\documentclass{amsart}

\usepackage[a4paper,hmargin=3.5cm,vmargin=4cm]{geometry}
\usepackage{amsfonts,amssymb,amscd,amstext}
\usepackage{graphicx}
\usepackage[dvips]{epsfig}
\usepackage{mathpazo}

\let\geq=\geqslant

\let\leq=\leqslant
\let\ptl=\partial
\let\Sg=\Sigma
\let\sg=\sigma
\let\eps=\varepsilon
\let\Om=\Omega

\newcommand{\subeq}{\subseteq}
\newcommand{\sub}{\subset}
\newcommand{\rr}{\mathbb{R}}
\newcommand{\rrn}{\mathbb{R}^{n+1}}

\newcommand{\nn}{\mathbb{N}}
\newcommand{\escpr}[1]{\left< #1\right>}

\DeclareMathOperator{\vol}{vol}
\DeclareMathOperator{\divv}{div}

\newtheorem{theorem}{Theorem}[section]
\newtheorem{proposition}[theorem]{Proposition}
\newtheorem{lemma}[theorem]{Lemma}
\newtheorem{corollary}[theorem]{Corollary}

\newtheorem{example}[theorem]{Example}

\theoremstyle{remark}

\newtheorem{remark}[theorem]{Remark}

\newtheorem{conjecture}[theorem]{\bf{Conjecture}}

\numberwithin{equation}{section}

\begin{document}

\title[The isoperimetric problem in $\rrn$ with density]{On the
isoperimetric problem in Euclidean space with density}

\author[C.~Rosales]{C\'esar Rosales}
\address{Departamento de
Geometr\'{\i}a y Topolog\'{\i}a \\ Facultad de Ciencias \\
Universidad de Granada \\ E-18071 Granada, Spain}
\email{crosales@ugr.es}

\author[A.~Ca\~nete]{Antonio Ca\~nete}
\address{Departamento de
Geometr\'{\i}a y Topolog\'{\i}a \\ Facultad de Ciencias \\
Universidad de Granada \\ E-18071 Granada, Spain}
\email{antonioc@ugr.es}

\author[V.~Bayle]{Vincent Bayle}
\address{Institute Fourier, BP 74,
38402 Saint Martin D'Heres Cedex, France}
\email{vbayle@ujf-grenoble.fr}

\author[F.~Morgan]{Frank Morgan}
\address{Department of Mathematics and Statistics \\
Williams College \\ Williamstown, MA 01267, U.~S.~A.}
\email{Frank.Morgan@williams.edu}

\thanks{First and second authors are partially supported by MCyT-Feder
research project MTM2004-01387, fourth author by the National Science
Foundation} \keywords{Manifolds with density, isoperimetric problem,
generalized mean curvature, stability, symmetrization}
\subjclass[2000]{49Q20, 53C17} \date{\today}

\begin{abstract}
We study the isoperimetric problem for Euclidean space endowed with a
continuous density.  In dimension one, we characterize isoperimetric
regions for a unimodal density.  In higher dimensions, we prove
existence results and we derive stability conditions, which lead to
the conjecture that for a radial log-convex density, balls about the
origin are isoperimetric regions.  Finally, we prove this conjecture
and the uniqueness of minimizers for the density $\exp (|x|^2)$ by
using symmetrization techniques.
\end{abstract}

\maketitle

\thispagestyle{empty}

\section{Introduction}
\label{sec:intro}

The \emph{isoperimetric problem} inside a Riemannian manifold seeks
regions of least perimeter enclosing a fixed amount of volume.  This
problem can also be studied in the more general setting of a
\emph{manifold with density}, where a given continuous positive
function on the manifold is used to weight the Riemannian volume and
boundary area.  Such a density is not equivalent to scaling the metric
conformally by a factor $\lambda$, since in that case volume and
perimeter would scale by different powers of $\lambda$.

We shall consider the particular case of Euclidean space with a
density $f=e^\psi$.  For any Borel set $\Om$ in $\rrn$, the (weighted)
\emph{volume} or \emph{measure} of $\Om$, and the (weighted)
\emph{perimeter relative to an open set} $U$ in $\rrn$ are given by
\[
\vol(\Om)=\int_{\Om}f\,dv,\qquad P(\Om,U)=\int_{\ptl\Om\cap U}f\,da,
\]
where $dv$ and $da$ are elements of Euclidean volume and area, in
general provided by Lebes\-gue measure and $n$-dimensional Hausdorff
measure on $\rrn$.  Let $P(\Om)=P(\Om,\rrn)$.

Much of the information of the isoperimetric problem is contained in
the \emph{isoperimetric profile}, which is the function
$I_{f}:(0,\vol(\rrn))\to\rr$ given by
\[
I_{f}(V)=\inf\,\{P(\Om): \Om \text{ is a smooth open set with
}\vol(\Om)=V\}.
\]
An \emph{isoperimetric region} --or simply a \emph{minimizer}-- of
volume $V$ is an open set $\Om$ such that $\vol(\Om)=V$ and
$P(\Om)=I_{f}(V)$.

In the last years the study of isoperimetric problems in manifolds
with density has increased.  One of the first and most interesting
examples, with applications in probability and statistics, is the
\emph{Gaussian density} $\exp(-\pi |x|^2)$.  About $1975$
C.~Borell~\cite{bor} and V.~N.~Sudakov and B.~S.~Tirel'son~\cite{st}
independently proved that half-spaces minimize perimeter under a
volume constraint for this density.  In $1982$ A. Ehrhard
\cite{Ehrhard} gave a new proof of the isoperimetric property of
half-spaces by adapting symmetrization techniques to the Gaussian
context.  His proof can be simplified by following the generalization
of Steiner symmetrization to product measures given by A~Ros
\cite{ros}.  More recently, S.~Bobkov and C.~Houdr{\'e} \cite{bh}
considered ``unimodal densities" with finite total measure in the real
line.  These authors explicitly computed the isoperimetric profile for
such densities and found some of the isoperimetric solutions.
M.~Gromov \cite{gromov} studied manifolds with density as ``mm
spaces'' and mentioned the natural generalization of mean curvature
obtained by the first variation of weighted area.  V.~Bayle \cite{ba}
proved generalizations of the L\'evy-Gromov isoperimetric inequality
and other geometric comparisons depending on a lower bound on the
generalized Ricci curvature of the manifold.  For recent advances on
manifolds with density we refer the reader to \cite{m4}, \cite{ba} and
references therein.

In this paper we first prove some existence results of isoperimetric
regions for densities in Euclidean space with infinite total measure
(Theorems~\ref{prop:exist1} and \ref{th:exist3}) and recall regularity
properties of minimizers (Theorem \ref{th:reg}).  In
Section~\ref{sec:variational} we use a variational approach to
characterize stability of balls centered at the origin for radial
densities (Theorem~\ref{th:stableballs}).  This result leads us to
Conjecture~\ref{conj:conjecture}: \emph{for radial log-convex
densities in $\rrn$, balls about the origin provide minimizers of any
given volume}.  We will prove this conjecture in the one-dimensional
case (Corollary~\ref{cor:logconvsym}) and for the radial density
$\exp(|x|^2)$ in any dimension (Theorem~\ref{th:main}).

In Section~\ref{sec:dim1} we completely describe isoperimetric regions
in the real line endowed with unimodal densities
(Theorems~\ref{th:fnodnoi} and~\ref{th:fnoinod}).  We use comparison
arguments that provide at the same time existence and uniqueness of
minimizers.  As interesting consequences we solve the isoperimetric
problem for log-concave and log-convex densities in the real line
(Corollaries~\ref{th:logconcave} and~\ref{th:cpletedescriplogconv}),
improving previous results by S. Bobkov and C. Houdr\'e~\cite{bh}.  We
also treat the isoperimetric problem and the free boundary problem for
the closed half-line $[0,+\infty)$ and for compact intervals.

In Section~\ref{sec:logconvdensity} we establish, in arbitrary
dimension, the isoperimetric property of round balls about the origin
for the density $\exp(c |x|^2)$, $c>0$ (Theorem~\ref{th:main}).  A
remarkable difference with respect to the Gaussian measure is that the
density $\exp(c |x|^2)$ for $c>0$ has infinite total volume and hence
the existence of minimizers is a non-trivial question.  The proof of
Theorem \ref{th:main} goes as follows.  First, we apply our previous
results in Section \ref{sec:prelimi} to ensure existence of
isoperimetric regions of any given volume.  Second, we use the
description of minimizers for log-convex densities on the real line
(Corollary~\ref{cor:logconvsym}) and the symmetrization in spaces with
product measures given by A.~Ros \cite{ros}, to construct a
counterpart to Steiner symmetrization for the density $\exp(c|x|^2)$.
Then we use this symmetrization in axis directions as employed by
L.~Bieberbach \cite{bieberbach} to produce centrally symmetric
minimizers with connected boundary.  Finally we conclude by Hsiang
symmetrization \cite{hsiang} that such a minimizer must be a round
ball about the origin.  As a corollary of Theorem \ref{th:main} we
deduce an eigenvalues comparison theorem for the density $\exp(c
|x|^2)$, $c\geq 0$, generalizing the Faber-Krahn Inequality.

Usually the uniqueness of isoperimetric regions is difficult to prove.
In the case of the Gaussian density, the complete characterization of
equality cases in the isoperimetric inequality is due to E.~A.~Carlen
and C.~Kerce \cite{ck}, who proved that any perimeter minimizer for
fixed volume is, up to a set of measure zero, a half-space.  They
obtained this result as consequence of the discussion of equality in a
more general functional inequality due to S.~Bobkov.  Previous
uniqueness results in the Gaussian setting involving a Brunn-Minkowski
type inequality were given by A.~Ehrhard \cite{Ehrhard2}.  In
Theorem~\ref{th:main} we also show the uniqueness of round balls
centered at the origin as minimizers for the density $\exp(c|x|^2)$,
$c>0$.  Since round balls appear as the result of finitely many
symmetrizations, it suffices to see that if an axis symmetrization of
a minimizer produces a ball, then the minimizer is a ball.  We deduce
this fact by standard arguments \cite[Lemma III.2.3]{chavel}.

An interesting consequence of our characterization of stable balls in
Theorem~\ref{th:stableballs} observed by K.~Brakke is that any round
ball about the origin is unstable in $\rrn$ endowed with a radial,
strictly log-concave density.  This fact, together with the
isoperimetric property of half-spaces in the Gaussian space and our
Corollary~\ref{cor:logconcave}, where we proved the isoperimetric
property of half-lines for densities on the real line, might suggest
that half-spaces are isoperimetric regions for any radial, log-concave
density on $\rrn$.  In Corollary~\ref{prop:counterexample} we give an
example showing that this is not true in general.  We believe that
this is a motivation to study in more detail the isoperimetric problem
for these kind of densities where unexpected shapes appear.

After circulating this manuscript we heard from Franck Barthe and
Michel Ledoux that the isoperimetric property of round balls about the
origin for the density $\exp{(c|x|^2)}$, $c>0$, was previously proved
by C.~Borell \cite[Theorem 4.1]{bor2}.  Borell's proof uses a
Brunn-Minkowski inequality and does not yield uniqueness of
minimizers.

\noindent\textbf{Acknowledgements.} This work began during Morgan's
lecturers on ``Geometric Measure Theory and Isoperimetric Problems" at
the 2004 Summer School on Minimal Surfaces and Variational Problems
held in the Institut de Math{\'e}matiques de Jussieu in Paris.  We
would like to thank the organizers: Pascal Romon, Marc Soret, Rabah
Souam, Eric Toubiana, Fr\'ed\'eric H\'elein, David Hoffman, Antonio
Ros and Harold Rosenberg.  We also thank Franck Barthe, Christer
Borell and Michel Ledoux for bringing previous results on the
isoperimetric problem for the density $\exp{(|x|^2)}$ to our
attention.

\section{Existence and regularity results}
\label{sec:prelimi}

In this section we firstly deal existence of isoperimetric regions in
Euclidean space with density.  In general, for a Riemannian manifold
with density, standard compactness arguments of Geometric Measure
Theory (see \cite[27.3 and 31.2]{simon} or \cite[5.5 and 9.1]{m1}, and
\cite[4.1]{m3} or \cite[Thm.  2.1]{ritros}) can be applied in order to
provide isoperimetric regions, except that there can be loss of volume
at infinity (by regularity Theorem~\ref{th:reg}, these are open sets
with nice boundaries).  In particular, if the total measure is finite,
isoperimetric regions of any prescribed volume exist.  We will prove
some existence results for densities with infinite total volume.  We
begin with the following lemma:

\begin{lemma}
\label{lem:exist2}
Let $f$ be a positive, nondecreasing function on
$[0,+\infty)$ satisfying $f(r)\to +\infty$ when $r\to +\infty$.  If
the function $\psi=\log(f)$ satisfies
\[
\psi(r)\leq C\left (\frac{n+1}{n}-\eps\right )^{r/2}
\]
for some $n\in\nn$, $C>0$ and $\eps\in (0,1)$, then the sequence
\[
\zeta(m)=\frac{f(m)}{f(m+2)^{n/(n+1)}}
\]
tends to infinity.

Conversely, if $\{\zeta(m)\}\to +\infty$, then there is $r_0>0$ and
$C>0$ such that
\[
\psi(r)\leq C\left(\frac{n+1}{n}\right)^{r/2},\qquad r\geq r_0.
\]
\end{lemma}

\begin{proof}
We prove the first part of the statement by contradiction.  Suppose
that $\{\zeta(m)\}_{m\in\nn}$ does not tend to infinity.  We can
assume, by passing to a subsequence if necessary, that there is
$K>0$ such that
\[
f(m)\leq K\,f(m+2)^{n/(n+1)},\qquad m\in\nn,
\]
and therefore
\[
\psi(m+2)\geq\frac{n+1}{n}\,(\psi(m)-\log(K)),\qquad m\in\nn.
\]
On the other hand, as $\{\psi(m)\}\to+\infty$ and $\eps>0$, we can
find $m_0\in\nn$ such that
\begin{align*}
\psi(m+2)&\geq\left
(\frac{n+1}{n}-\frac{\eps}{2}\right)\,\psi(m),\qquad m\geq m_0,
\\
\psi(m_{0}+2k)&\geq\left
(\frac{n+1}{n}-\frac{\eps}{2}\right)^k\,\psi(m_{0}).
\end{align*}
Now take $r\geq m_0+2$ and $k\in\nn$ such that $m_{0}+2k\leq r\leq
m_{0}+2(k+1)$.  By using that $\psi$ is nondecreasing we have
\begin{align*}
\psi(r)\geq\psi(m_{0}+2k)\geq\left(\frac{n+1}{n}-\frac{\eps}{2}\right)^{k}
\,\psi(m_{0})\geq\left(\frac{n+1}{n}-\frac{\eps}{2}\right)^{r/2-m_{0}/2-1}\,
\psi(m_{0}).
\end{align*}
Hence, for $r\gg m_0+2$ we deduce
\[
\psi(r)>C\,\left(\frac{n+1}{n}-\eps\right)^{r/2},
\]
and we get a contradiction.

Conversely, suppose that $\{\zeta(m)\}\to+\infty$.  Then, we can find
$m_0\geq 2$ such that $\psi(m_0)>0$ and $\zeta(m)\geq 1$ for $m\geq
m_0$.  As a consequence
\begin{align*}
\psi(m+2)&\leq\left(\frac{n+1}{n}\right)\psi(m),\qquad m\geq m_0,
\\
\psi(m_0+2k)&\leq\left(\frac{n+1}{n}\right)^k\,\psi(m_0).
\end{align*}
Finally, for $r\geq m_0$ there is $k\in\nn$ such that $m_0+2(k-1)\leq
r\leq m_0+2k$.  Hence
\[
\psi(r)\leq\psi(m_0+2k)\leq\left(\frac{n+1}{n}\right)^k\,\psi(m_0)
\leq\left(\frac{n+1}{n}\right)^{r/2}\,\psi(m_0).
\]
\end{proof}

Now, we can prove our first existence result.

\begin{theorem}
\label{prop:exist1}
Let $f=e^\psi$ be a density on $\rrn$ such that $f(x)\to +\infty$ when
$|x|\to +\infty$.  Suppose that one of the following conditions holds:
\begin{itemize}
\item[(i)] The sequence defined by
\[
\zeta(m)=\frac{\min\,\{f(x):m\leq |x|\leq m+2\}}
{\max\,\{f(x)^{n/(n+1)}:m\leq |x|\leq m+2\}}
\]
tends to infinity.
\item[(ii)] The density is radial, nondecreasing in $|x|$ and satisfies
\[
\psi(x)\leq C\,\left(\frac{n+1}{n}-\eps\right)^{|x|/2},
\]
for some constants $C>0$ and $\eps\in (0,1)$.
\end{itemize}
Then, minimizers of any given volume exist for this density and they are bounded subsets of
$\rrn$.
\end{theorem}

\begin{remark}
The proof of the statement shows that it suffices to suppose that the
ratio $\min f(x)/\max f(x)^{n/(n+1)}$ on lattice cubes
goes to infinity.
\end{remark}

\begin{proof}
By Lemma~\ref{lem:exist2} we can assume that (i) holds.  Denote by $v(\Om)$ and $a(\ptl\Om)$
the Euclidean volume and boundary area of a set $\Om$.  Partition $\rrn$ into lattice open
cubes of diameter equal to $1$ and Euclidean volume $v_{0}$.  There is an isoperimetric
constant $\alpha>0$ such that any set $\Om$ inside a cube $C$ as above with $v(\Om)\leq
v_{0}/2$ satisfies
\[
a(\ptl\Om\cap C)\geq\alpha\,v(\Om)^{n/(n+1)}.
\]
On the other hand, there is $m=m(C)\in\nn$ such that the cube $C$ is
contained in the annulus $\{m\leq |x|\leq m+2\}$.  Thus, the
definition of weighted volume and perimeter, together with the
definition of $\zeta(m)$, implies the inequality
\begin{equation}
\label{eq:chok1}
P(\Om,C)\geq\alpha\,\zeta(m)\vol(\Om)^{n/(n+1)},
\end{equation}
for any $\Om\sub C\sub\{m\leq |x|\leq m+2\}$ with $v(\Om)\leq
v_{0}/2$.

Fix $V>0$, and consider a sequence of smooth open sets of volume $V$ with perimeters
approaching $I_f(V)$ and bounded from above by $I_{f}(V)+1$.  By using the Compactness Theorem
\cite[9.1]{m2} we may assume that this sequence converges.  Fix $\eps>0$.  By hypothesis, there
is $m_{0}\in\nn$ such that $\zeta(m)\geq (1/\eps)^{n/(n+1)}$ for any $m\geq m_{0}$.  On the
other hand, as $f(x)\to +\infty$ when $|x|\to +\infty$, we can suppose that $v(\Om)\leq
v_{0}/2$ whenever $\Om\sub C\sub\{|x|\geq m_{0}\}$.  In particular, we can apply
\eqref{eq:chok1} to such an $\Om$, so that we obtain
\[
\vol(\Om)\leq\left (\frac{P(\Om,C)}{\alpha\,\zeta(m)}
\right)^{(n+1)/n}\leq\eps\,\alpha'P(\Om,C)^{(n+1)/n}.
\]
By summing the previous inequality over the collection $\mathcal{C}_m$ of all cubes $C$
contained in $\{|x|\geq m\}$ we deduce that for any set $\Om$ of the given minimizing sequence
and any $m\geq m_{0}$
\begin{align*}
\vol\bigg(\Om\cap(\bigcup_{C\in\,\mathcal{C}_m}C)\bigg)&\leq\eps\,\alpha'\left
(\sum_{C\in\,\mathcal{C}_m} P(\Om,C) \right )^{(n+1)/n}
\\
&\leq\eps\,\alpha'P(\Om)^{(n+1)/n}
\leq\eps\,\alpha'\,(I_{f}(V)+1)^{(n+1)/n}.
\end{align*}
Hence, there is no loss of volume at infinity and the limit of our
sequence is an isoperimetric region of volume $V$.

To prove that any minimizer $\Om$ is a bounded subset of $\rrn$, we
can proceed as in \cite[Lemma 13.6]{m1}.  Consider any large cube
$C_r=[-r,r]^{n+1}$ about the origin and partition almost all its
complement into congruent open cubes of diameter at most $1$.  Denote
$V(r)=\vol(\Om-C_r)$ and $P(r)=P(\Om,\rrn-C_r)$.  As above,
we have
\begin{equation}
\label{eq:bounded1} V(r)\leq\eps\,\alpha' P(r)^{(n+1)/n},\qquad r\gg 0.
\end{equation}
On the other hand, there is a constant $H>0$ depending on $\ptl\Om$
such that small volume adjustments may be accomplished inside $C_r$ at
a cost
\[
|\Delta P|\leq H\,|\Delta V|.
\]
Thus replacing $\Om-C_r$ costs at most $H|\Delta V|+|V'(r)|$ (due to the slice of $\ptl C_r$)
for almost all large $r$. By using that $\Om$ is a minimizer, we get
\begin{equation}
\label{eq:bounded2} P(r)\leq H\,V(r)+|V'(r)|,\qquad\text{ for almost all }\, r\gg 0.
\end{equation}
Since $V(r)$ is nonincreasing and tends to $0$ when $r\to +\infty$, combining inequalities
\eqref{eq:bounded1} and \eqref{eq:bounded2} yields for some $c>0$,
\[
c\,V(r)^{n/(n+1)}\leq -V'(r),\qquad\text{ for almost all }\, r\gg 0.
\]
If we suppose that $\Om$ is unbounded, then $V(r)\neq 0$ and
\[
(n+1)\,(V^{1/(n+1)})'=V^{-n/(n+1)}\,V'\leq -c<0,
\]
for almost all large $r>0$, a contradiction since $V$ is positive and nonincreasing.
\end{proof}

Our next existence result is an improvement of Theorem~\ref{prop:exist1} in dimension two. We
need the following lemma:

\begin{lemma}
\label{lem:desisoper}
Let $f$ be a planar radial density nondecreasing on $[r_{0},+\infty)$.
Then, for any smooth, open set $\Om\sub\rr^2$ contained in $\{|x|\geq
r_{0}\}$ and such that $P(\Om)<2\pi r_{0} f(r_{0})$, we have the
isoperimetric inequality
\[
P(\Om)^2\geq 2\,f(r_{0})\vol(\Om).
\]
\end{lemma}

\begin{proof}
First, we can assume that $\Om$ is connected.  Moreover, the hypothesis on the perimeter
implies that $\Om$ is bounded and the closure of $\Om$ cannot contain a circle about the
origin.  Let $r_{1}$ and $r_{2}$ be the minimum and maximum distance from $\overline{\Om}$ to
the origin, respectively.  The intersection $\Om_{t}$ of $\Om$ with the circle of radius $t\in
(r_{1},r_{2})$ has Euclidean length strictly less than $2\pi t$, and the boundary $\ptl\Om_{t}$
has at least two points.  Therefore, the coarea formula gives us
\begin{equation}
\label{eq:dios1}
P(\Om)\geq\int_{r_{1}}^{r_{2}}f(t)\,\text{card} (\ptl\Om_{t})\,dt\geq
2f(r_{0})\,(r_{2}-r_{1}),
\end{equation}
where we have used that the density is nondecreasing on
$[r_{0},+\infty)$.

On the other hand, we consider the map $F:(r_{1},r_{2})\times\ptl\Om\to\rr^2$ given by
$F(t,x)=t x/|x|$. It is clear that $\Om\subeq F(A)$, where $A$ is the open set of the pairs
$(t,x)$ where $t<|x|$ and $f(tx/|x|)<f(x)$. For any $(t,x)\in A$ the Jacobian of $F$ is
strictly less than $1$. Thus, the definition of $A$, together with the coarea formula and
Fubini's theorem implies
\begin{equation}
\label{eq:dios2}
\vol(\Om)\leq\vol(F(A))\leq\int_{A}f\left(\frac{tx}{|x|}\right)
\,d(t,x)\leq\int_{A}f(x)\,d(t,x)=(r_{2}-r_{1})\,P(\Om).
\end{equation}
Multiplying the estimates \eqref{eq:dios1} and \eqref{eq:dios2} we
obtain the desired inequality.
\end{proof}

\begin{remark}
We do not see how to generalize the previous lemma to $\rrn$. The
analog of inequality \eqref{eq:dios2} holds, but the estimation on the
Euclidean boundary area $a(\ptl\Om_{t})$ leading to \eqref{eq:dios1} becomes
\[
a(\ptl\Om_{t})\geq\alpha\,f(t)^{1/n}\,v(\Om_{t})^{(n-1)/n}\geq
C\,v(\Om_{t}),
\]
where $v(\Om_{t})$ denotes Euclidean volume. The last inequality can be
integrated to deduce $P(\Om)\geq C\,\vol(\Om)$, which is inadequate
to obtain an analog of Lemma \ref{lem:desisoper}.
\end{remark}

\begin{theorem}
\label{th:exist3}
Consider the plane endowed with a nondecreasing, radial density $f$ such
that $f(x)\to +\infty$ when $|x|\to +\infty$. Then, there are
minimizers for this density of any given volume.
\end{theorem}

\begin{proof}
Consider a sequence of smooth open sets of volume $V>0$ with perimeters
approaching $I_f(V)$.  Applying the Compactness Theorem \cite[9.1]{m2}
we can assume that this sequence converges.  We can also suppose that any set
$\Om$ of this sequence satisfies $P(\Om)\leq I_f(V)+1$.  Moreover, as
the density tends to $+\infty$, there is $m_{0}\in\nn$ such that
$I_f(V)+1<2\pi m f(m)$ for any $m\geq m_{0}$.  In particular, we can
apply Lemma \ref{lem:desisoper} to the union $\Om'$ of all connected
components of $\Om$ inside $\{|x|\geq m\}$.  Hence, we get
\[
\vol(\Om')\leq\frac{P(\Om')^2}{2f(m)}\leq\frac{(I_{f}(V)+1)^2}{2f(m)}
,\qquad m\geq m_{0}.
\]
As $\lim_{m\to +\infty} f(m)=+\infty$ we conclude that there is no
loss of volume at infinity and the limit of our sequence solves the
isoperimetric problem for volume $V$.
\end{proof}

\begin{example}
\emph{We illustrate here that Theorem \ref{th:exist3} need not hold if
we do not require the density to be nondecreasing.  Consider in $\rrn$
$(n>1)$ the density $f(x)=1+|x|^2$.  Now, introduce bumps into the
graph of $f$ such that any volume $V_{k}$ corresponding to a positive
rational can be enclosed with perimeter $1/k$.  Then, for any given
volume we may find a sequence of sets enclosing this volume and with
arbitrarily small perimeter, which implies that isoperimetric regions
do not exist.}
\end{example}

We finish this section by recalling regularity properties of the boundary of a minimizer in
Euclidean space with density.  The result is also valid for any smooth Riemannian manifold with
density.

\begin{theorem}[{\cite[3.10]{m2}}]
\label{th:reg} Consider a smooth density on $\rrn$.  If $\Om$ is a minimizer, then the boundary
$\Sg=\ptl\Om$ is a real-analytic embedded hypersurface, up to a closed set of singularities
with Euclidean Hausdorff dimension less than or equal to $n-7$.
\end{theorem}

\section{Variational formulae. Stable balls for radial densities}
\label{sec:variational}

In this section we use a variational approach to derive some
properties of sets minimizing perimeter up to second order
for variations preserving volume.

Let $f=e^\psi$ be a smooth density on $\rrn$.  Denote by $\Om\sub\rrn$ a smooth open set with
boundary $\Sg$ and inward unit normal vector $N$.  We consider a one-parameter variation
$\{\phi_t\}_{|t|<\eps}:\rr^{n+1}\rightarrow\rr^{n+1}$ with associated infinitesimal vector
field $X=d\phi_t/dt$ with normal component $u=\escpr{X,N}$.  Let $\Om_t=\phi_t(\Om)$ and
$\Sigma_t=\phi_t(\Sigma)$.  The volume and perimeter functions of the variation are
$V(t)=\vol(\Om_{t})$ and $P(t)=P(\Om_{t})$, respectively.  The first variation of volume and
perimeter are computed in \cite[Chapter 3]{ba}.  We include here a proof for the sake of
completeness.

\begin{lemma}
\label{lem:firstvariation}
The first variation of volume and perimeter of a smooth region $\Om$
with boundary $\Sg$ in $\rrn$ endowed with smooth density $f=e^\psi$ for a
flow with initial normal velocity $u$ are given by
\[
V'(0)=-\int_\Sigma fu\,dv,\qquad
P'(0)= -\int_\Sigma (nH-\escpr{\nabla\psi,N})\,fu\,da,
\]
where $H$ is the Euclidean mean curvature of $\Sigma$ with respect to
$N$ (that is, the arithmetic mean of the principal curvatures of
$\Sg$) and $\nabla\psi$ is the Euclidean gradient of $\psi$.
\end{lemma}

\begin{proof}
Denote by $\divv X$ (resp.  $\divv_{\Sg} X$) the divergence of $X$ in
$\rrn$ (resp.  relative to $\Sg$).  Let $X(f)=\escpr{\nabla f,X}$.  We
have
\begin{align*}
V'(0)&=\int_{\Om}X(f)\,dv+\int_{\Om}
f\,\frac{d}{dt}\bigg |_{t=0} (dv_{t})
\\
&=\int_\Om(\escpr{\nabla f,X}+f\,\divv{X})\,dv=
\int_\Om\divv(fX)\,dv=-\int_\Sg fu\,da.
\end{align*}
In the second equality we have used that
$(d/dt)|_{t=0}\,(dv_{t})=(\divv X)\,dv$, see \cite[\S 16]{simon}.  In
the last one, we have applied the Gauss-Green theorem.  For perimeter
we have
\begin{align*}
P'(0)&=\int_\Sg X(f)\,da+\int_{\Sg}f\,\frac{d}{dt}\bigg |_{t=0}
(da_{t})
\\
&=\int_{\Sg}(\escpr{\nabla f,X}+f\,\divv_\Sg{X})\,da=\int_\Sg(\escpr{\nabla
f,uN}+\divv_\Sg(fX))\,da
\\
&=\int_\Sg fu\escpr{\nabla\psi,N}\,da-\int_\Sg nHfu\,da =-\int_\Sg
(nH-\escpr{\nabla\psi,N})\,fu\,da.
\end{align*}
To obtain the fourth equality we have used that the integral over
$\Sg$ of the divergence of the tangent part of $fX$ vanishes by virtue
of the divergence theorem.
\end{proof}

We define, as in \cite[Chapter 3]{ba}, the (generalized) \emph{mean
curvature} of $\Sigma$ with respect to $N$ as the function
\begin{equation}
\label{eq:mcdef}
H_\psi=nH-\escpr{\nabla\psi,N},
\end{equation}
so that the first variation of perimeter can be written as
\[
P'(0)=-\int_\Sigma H_\psi fu\,da.
\]

We say that a given variation $\{\phi_{t}\}_t$ \emph{preserves volume} if $V(t)$ is constant
for any small $t$.  We say that $\Om$ is \emph{stationary} if $P'(0)=0$ for any
volume-preserving variation. It is clear that any isoperimetric region is also stationary.  The
following characterization of stationary sets is similar to the one established by
J.~L.~Barbosa and M.~do Carmo \cite[Proposition 2.7]{bdc} for the case of $\rrn$ with the
standard density $f\equiv 1$.  The proof is based on Lemma~\ref{lem:firstvariation} and on the
fact that any function $u$ orthogonal to $f$ in $L^2(\Sg)$ is the normal component of a vector
field associated to a volume-preserving variation of $\Om$, see \cite[Lemma 2.2]{bdc}.

\begin{proposition}
\label{prop:stationary} Consider a smooth density $f=e^\psi$ on $\rrn$. Then, for a smooth open
set $\Om$, the following conditions are equivalent:
\begin{itemize}
\item[(i)] $\Om$ is stationary.
\item[(ii)] $\Sg=\ptl\Om$ has $($generalized$)$ constant mean curvature
$H_{0}$.
\item[(iii)] There is a constant $H_0$ such that $(P-H_0V)'(0)=0$ for any
variation of $\Om$.
\end{itemize}
\end{proposition}

\begin{example}
\label{ex:real} \emph{Let $f=e^\psi$ be a smooth density defined on the real line. Then, it is
easy to show that a bounded interval $(a,b)$ is stationary if and only if
$\psi'(a)=-\psi'(b)$.}
\end{example}

Now, we introduce some examples of hypersurfaces with constant mean
curvature in $\rrn$ with a radial density.

\begin{example}
\label{re:mcspheres} \emph{Suppose $f=e^\psi$, where $\psi(x)=\delta(|x|)$ for any $x\in\rrn$.
The mean curvature $H_{\psi}$ of a hypersurface $\Sg$ with respect to a unit normal vector $N$
is given by
\[
H_\psi(p)=nH(p)-\frac{\delta'(r)}{r}\,\escpr{p,N(p)},\quad r=|p|,
\]
where $H$ is the Euclidean mean curvature of $\Sg$ with respect to $N$. In particular, if $\Sg$
is a sphere of radius $r>0$, then it has constant mean curvature if and only if $\Sg$ is
centered at the origin.  In this case, $H_\psi=n/r+\delta'(r)$ with respect to the inner normal
vector.}

\emph{On the other hand, if $\Sg$ is the hyperplane defined by
$\{x\in\rr^{n+1}:\escpr{x,u}=c\}$, where $|u|=1$, then the mean curvature of $\Sg$ with respect
to $N=-u$ is
\begin{equation}
\label{eq:mc2}
H_\psi(p)=-c\,\,\frac{\delta'(r)}{r},\quad r=|p|.
\end{equation}
It follows that any hyperplane passing through the origin is a minimal hypersurface of $\rrn$
with a radial density.  In general, we cannot expect that any hyperplane has constant mean
curvature for a radial density.  In fact, a straightforward analysis of equation \eqref{eq:mc2}
leads us to the following:}
\end{example}

\begin{lemma}
\label{lem:puc} Let $f=e^\psi$ be a smooth radial density on $\rrn$.  Suppose that there is a
hyperplane $\Sg$ which does not contain the origin and has constant mean curvature $H_{\psi}$.
Then, there are constants $a,b\in\rr$ and $r_{0}>0$, such that
\[
\psi(x)=e^{\,a |x|^2+b},\quad\text{ whenever }|x|\geq r_0.
\]
\end{lemma}

Now, we compute the second variation formula of the functional $P-H_{\psi}V$ for any variation
of a stationary set.

\begin{proposition}[{\cite[Section 3.4.6]{ba}}]
Consider a stationary open set $\Om$ in $\rrn$ endowed with a smooth
density $f=e^\psi$.  Let $N$ be the inward unit normal vector to
$\Sg=\ptl\Om$, and $H_{\psi}$ the constant mean curvature of $\Sg$
with respect to $N$.  Consider a variation of $\Om$ with associated
vector field $X=~uN$ on $\Sg$.  Then, we have
\begin{equation}
\label{eq:second}
(P-H_{\psi}\,V)''(0)=Q_\psi(u,u):=\int_\Sigma
f\,(|\nabla_\Sigma u|^2-|\sigma|^2u^2)\,da +\int_\Sigma
fu^2\,(\nabla^2\psi)\,(N,N)\,da,
\end{equation}
where $\nabla_\Sigma u$ is the gradient of $u$ relative to $\Sigma$,
$|\sigma|^2$ is the squared sum of the principal curvatures of
$\Sigma$, and $\nabla^2\psi$ is the Euclidean Hessian of $\psi$.
\end{proposition}

\begin{proof}
The first variation formula for volume and perimeter gives us
\[
(P-H_\psi V)'(t)=-\int_{\Sg_t}(H_\psi)_t\,
fu_{t}\,da_t+H_\psi\int_{\Sg_t}fu_{t}\,da_t,
\]
where $(H_{\psi})_{t}$ is the mean curvature of $\Sg_{t}$. Hence
\begin{equation}
\label{eq:tocado}
(P-H_\psi V)''(0)=-\int_\Sg H_\psi'(0)\,fu\,da,
\end{equation}
so that we have to compute the derivative of the generalized mean
curvature along $\Sg_t$.  Denote by $D_{U}V$ the Levi-Civit\'a
connection on $\rrn$.  By \eqref{eq:mcdef}, we get
\begin{align*}
H_\psi'(0)&=n H'(0)-\escpr{D_{X}\nabla\psi,N}-
\escpr{\nabla\psi,D_{X} N}
\\
&=n H'(0)-u\,(\nabla^2\psi)(N,N)+\escpr{\nabla\psi,\nabla_{\Sg}u},
\end{align*}
where in the last equality we have used that $D_{X}N=-\nabla_{\Sg}u$.
On the other hand, it is well known \cite{rosenberg} that
\begin{equation}
\label{eq:loko10}
nH'(0)=\Delta_\Sg u+|\sg|^2u,
\end{equation}
where $\Delta_\Sg$ is the Laplacian relative to $\Sg$. Thus, we have
obtained
\[
H_\psi'(0)=\Delta_\Sg u+|\sg|^2u-u\,(\nabla^2\psi)(N,N)+
\escpr{\nabla_{\Sg}\psi,\nabla_{\Sg}u}.
\]
By substituting this information into \eqref{eq:tocado} we conclude that
\begin{align*}
(P-H_\psi V)''(0)=&-\int_\Sg fu\,(\Delta_\Sg u+|\sg|^2 u)\,da-\int_\Sg
fu\,\escpr{\nabla_\Sg\psi,\nabla_\Sg u}\,da
\\
&+\int_\Sg fu^2\,(\nabla^2\psi)\,(N,N)\,da.
\end{align*}
Finally, by using integration by parts, we deduce
\[
-\int_\Sg fu\,(\Delta_\Sg u+|\sg|^2 u)\,da-\int_\Sg fu\,\escpr{\nabla_\Sg\psi,\nabla_\Sg
u}\,da=\int_\Sigma f\,(|\nabla_\Sigma u|^2-|\sigma|^2u^2)\,da,
\]
and the result follows.
\end{proof}

\begin{remark}
In a smooth Riemannian manifold with density the second variation has an additional term
depending on the Ricci curvature of the manifold in the normal direction $N$. This term comes
from \eqref{eq:loko10} and it is given by
\[
-\int_{\Sg}\text{Ric}(N,N)\,fu^2\,da.
\]
\end{remark}

The expression \eqref{eq:second} defines a quadratic form on $C_0^\infty(\Sg)$ called the
\emph{index form} associated to $\Sg$. We say that a smooth open set $\Om$ is \emph{stable} if
it is stationary and $P''(0)\geq 0$ for any volume-preserving variation of $\Om$.  Stability
can be characterized in terms of the index form as in \cite[Proposition 2.10]{bdc}. More
precisely, we have the following:

\begin{lemma}
\label{le:stable}
Let $\Om$ be a smooth open set in $\rrn$ endowed with a smooth density
$f=e^\psi$.  Then, $\Om$ is stable if and only if it is stationary and
the index form \eqref{eq:second} of $\Sg=\ptl\Om$ satisfies
\[
Q_\psi(u,u)\geq 0 \ \text{ for any }u\in C^\infty_0(\Sigma)\text{
such that }\int_\Sigma fu\,da=0.
\]
\end{lemma}

Observe that the term in the index form containing $\nabla^2\psi$ indicates that the notion of
stability is more restrictive when the density $f=e^\psi$ is log-concave.  In fact, by
inserting in \eqref{eq:second} locally constant nowhere vanishing functions we easily deduce

\begin{corollary}
\label{lem:connectedness}
If $\Om$ is a smooth stable region in $\rrn$ with a smooth,
log-concave density, then the hypersurface $\Sg=\ptl\Om$ is connected
or totally geodesic.  Moreover, if the density is strictly
log-concave, then $\Sg$ is connected.
\end{corollary}

Our main result in this section characterizes the stability of
round balls about the origin for radial densities.

\begin{theorem}
\label{th:stableballs}
Consider a smooth density $f=e^\psi$ on $\rrn$
such that $\psi(x)=\delta(|x|)$.  Then, the round ball $B$ about the
origin of radius $r>0$ is stable if and only if $\delta''(r)\geq 0$.
\end{theorem}

\begin{proof}
We use Lemma \ref{le:stable}.  Denote by $\Sg$ the boundary of $B$,
and by $N$ the inward unit normal vector to $\Sg$.  Clearly the
density is constant on $\Sg$, so that a function $u\in C^\infty(\Sg)$
is orthogonal to $f$ in $L^2(\Sg)$ if and only if it has mean zero on
$\Sg$.  Moreover, $(\nabla^2\psi)(N,N)=\delta''(r)$ on $\Sg$.  As
consequence, the index form \eqref{eq:second} is given by
\[
Q_{\psi}(u,u)=f(r)\,\int_{\Sg}(|\nabla_\Sg
u|^2-|\sg|^2u^2)\,da+f(r)\,\delta''(r)\int_{\Sg}u^2\,da.
\]

Since Euclidean balls are stable regions in $\rrn$ with the standard
density $f\equiv 1$, the first integral is
nonnegative and vanishes for translations.  Consequently, if
$\delta''(r)\geq 0$, then $Q_\psi\geq 0$ and $B$ is stable.
Conversely, if $B$ is stable under infinitesimal translations, then
$\delta''(r)\geq 0$.
\end{proof}

As an immediate consequence of Theorem~\ref{th:stableballs} we obtain

\begin{corollary}
\label{cor:logconvexstable}
In Euclidean space endowed with a smooth,
radial, log-convex density, round balls centered at the origin are
stable regions.
\end{corollary}

The preceding corollary leads to the following conjecture inspired by
Ken Brakke at Jussieu:

\begin{conjecture}
\label{conj:conjecture}
\emph{In $\rrn$ with a smooth, radial, log-convex density, balls about the
origin provide isoperimetric regions of any given volume.}
\end{conjecture}

In Sections~\ref{sec:dim1} and \ref{sec:logconvdensity} we will prove
some special cases of this conjecture.  Another interesting
consequence of Theorem~\ref{th:stableballs} is the fact that for a
strictly log-concave density on $\rrn$, round balls about the origin
are unstable.  This allows us to prove the following:

\begin{corollary}
\label{prop:counterexample}
There are smooth, radial, log-concave
densities with finite volume in $\rr^2$ for which isoperimetric
regions are neither half-planes nor round balls.
\end{corollary}

\begin{proof}
Consider the density $f=e^\psi$, with $\psi(x)=-\sqrt{|x|^2+1}$.  The
total volume of this density is finite and hence minimizers of any
given volume exist, as was indicated at the beginning of
Section~\ref{sec:prelimi}.  The Hessian of $\psi$ at any
$(x,y)\in\rr^2$ is given by
\[
(\nabla^2\psi)_{(x,y)}(a,b)=\frac{-(bx-ay)^2-a^2-b^2}{(1+x^2+y^2)^{3/2}},
\]
and hence $f$ is strictly log-concave.  It follows by Example~\ref{re:mcspheres} and
Theorem~\ref{th:stableballs} that any round disk is unstable for this density.  On the other
hand, by taking into account Example~\ref{re:mcspheres} and Lemma~\ref{lem:puc}, we deduce that
only planes passing through the origin have constant mean curvature $H_\psi$.  As consequence,
a minimizer with measure different from the half volume of $\rr^2$ cannot be a disk nor a
half-plane.
\end{proof}

\section{Isoperimetry in the real line with density}
\label{sec:dim1}

In this section we study the isoperimetric problem in the real line with a \emph{unimodal
density}: a density which is increasing (or decreasing) on $(-\infty,x_0)$ and decreasing (or
increasing) on $(x_0,+\infty)$, for some $x_0\in (-\infty,+\infty]$.  S.~Bobkov and C.~Houdr\'e
\cite[Section 13]{bh} previously considered this setting under the further assumption of finite
total measure.  They computed the isoperimetric profile and gave some examples of isoperimetric
regions.  We provide here a simple, more general approach, which leads us to the complete
description of minimizers.

We begin by solving the isoperimetric problem for monotonic densities.
We recall that, for a function $f$, an end $E=\pm\infty$ has
\emph{finite measure} if $f$ is integrable in a neighborhood of $E$.

\begin{proposition}
\label{prop:fmonotona}
Let $f$ be a monotonic density on $\rr$ and denote by $E$ the end
where $f$ attains its infimum.  If $E$ has finite measure, then for
any given volume, a half-line containing $E$ is the unique isoperimetric
region.  If $E$ has infinite measure, then the isoperimetric profile
coincides with $2\,f(E)$, and it is approached or attained by a
bounded interval going off to $E$.
\end{proposition}

\begin{proof}
If $E$ has finite measure, then any candidate other than the half-line of the same measure has
at least two boundary points and hence greater perimeter since at least one of them is beyond
the half-line. If $E$ has infinite measure, then any open set enclosing a given volume has at
least two boundary points, so that the infimum perimeter is $2\,f(E)$, approached or attained
as asserted.
\end{proof}

\begin{example}
\emph{For the density $f(x)=e^x$ the end $E=-\infty$ has finite measure. Then the half-lines
$(-\infty,x)$ are the unique minimizers for fixed volume and the isoperimetric profile is given
by $I_{f} (V)=V$, for any $V>0$. For the density $f\equiv 1$ the profile is constant and
isoperimetric regions are bounded intervals.  Finally, the density $f(x)=e^x+1$ is one for
which the profile is constant while minimizers do not exist.}
\end{example}

We say that a function $f$ is \emph{increasing-decreasing} if there is $x_0\in\rr$ such that
$f$ is increasing on $(-\infty,x_0)$ and decreasing on $(x_0,+\infty)$, not necessarily
strictly.

\begin{theorem}
\label{th:fnodnoi} Let $f$ be an increasing-decreasing density on $\rr$. Then, if a minimizer
for given volume exists, it is a half-line, a bounded interval where $f$ attains its maximum or
equals its one-sided minimum, or the complement of one of these intervals. If it does not
exist, the infimum perimeter is approached by a bounded interval going off to $\pm\infty$.
\end{theorem}

\begin{proof}
Consider a smooth open set of the prescribed volume. If the closure contains a point $x_{0}$
where $f$ attains its maximum, then we can replace the given set with an interval containing
$x_{0}$. Otherwise, we can assume by Proposition \ref{prop:fmonotona} that the set consists of
one interval on one side, or an interval on each side of the maxima of $f$. Among intervals
containing $x_{0}$ there is one of least perimeter. Among intervals on one side, the infimum
perimeter is $2\min\{f(-\infty),f(+\infty)\}$ or the unique minimizer is a half-line. The
theorem follows.
\end{proof}

In the following examples we illustrate that the different possibilities in
Theorem~\ref{th:fnodnoi} may occur.

\begin{figure}[h]
\begin{center}
\includegraphics[angle=-90,width=8cm]{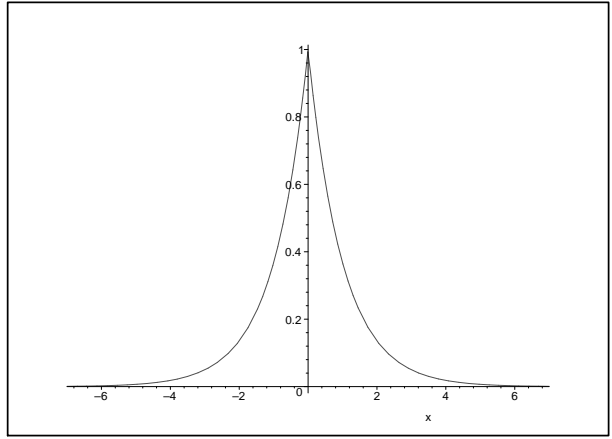}
\caption{\emph{For density $e^{-|x|}$ all types of minimizers occur.}}
\label{fig2}
\end{center}
\end{figure}

\begin{example}
\label{re:houseroof} \emph{Consider the density $f(x)=e^{-|x|}$, which has finite total measure
(see Fi\-gure~\ref{fig2}).  A straightforward computation shows that the isoperimetric
candidates of volume $V=1$ (half-lines, bounded intervals containing the origin, and
complements) have the same perimeter.  This illus\-trates that there is no uniqueness of
minimizers for $V=1$ since the different possibilities appear.  Moreover, though the density is
symmetric with respect to the origin, the bounded minimizers need not be symmetric.}
\end{example}

\begin{figure}[h]
\begin{center}
\includegraphics[angle=-90,width=8cm]{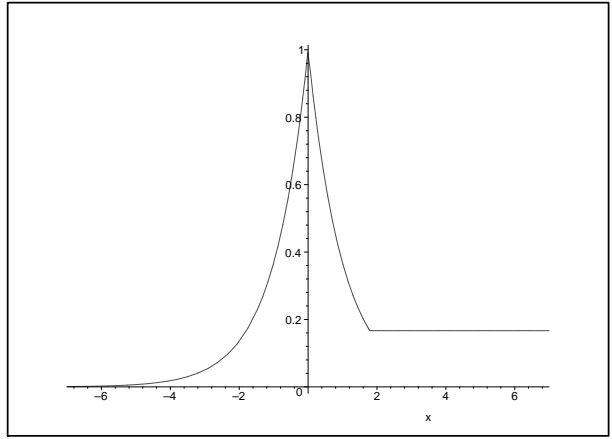}
\caption{\emph{A density for which minimizers are half-lines
and bounded intervals.}}
\label{fig3}
\end{center}
\end{figure}

\begin{example}
\emph{Consider the density given by $f(x)=e^{-|x|}$ for $x\leq\log(6)$ and $f(x)=1/6$ for
$x\geq\log(6)$.  The left end has finite measure while the right one has infinite measure (see
Figure~\ref{fig3}). Thus, only half-lines containing $-\infty$ and bounded intervals are
possible minimizers for a fixed measure.  For volume $V=1/3$, it can be shown that the
isoperimetric regions are the corresponding half-line containing $-\infty$ and any bounded
interval inside the half-line $[1/6,+\infty)$.  For a volume $V>1/3$ only bounded intervals
contained in $[1/6,+\infty)$ provide minimizers.  This illustrates that bounded minimizers need
not contain a point where the maximum of the density is achieved.}
\end{example}

\begin{figure}[h]
\begin{center}
\includegraphics[angle=-90,width=8cm]{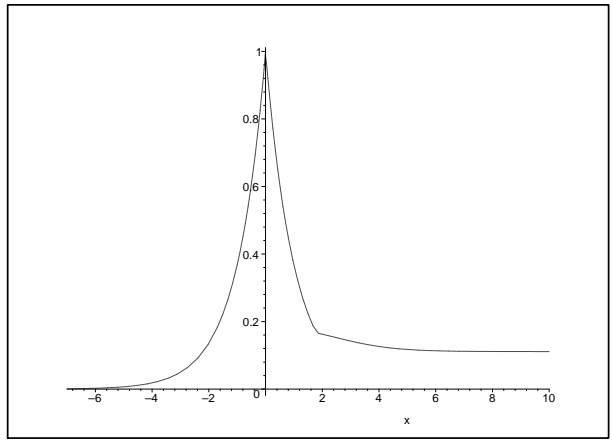}
\caption{\emph{A density for which minimizers do not always exist.}}
\label{fig4}
\end{center}
\end{figure}

\begin{example}
\emph{Consider the density given by
\[
f(x)=
\begin{cases}
e^{-|x|}\qquad\qquad\qquad x\leq\log(6),
\\
\frac{1}{9}+\frac{1}{x-\log(6)+18}\qquad\!\!\!  x\geq\log(6),
\end{cases}
\]} \!\!\emph{which is depicted in Figure~\ref{fig4}.  As in the
previous example, the ends $-\infty$ and $+\infty$ have finite and infinite measure,
respectively.  It is not difficult to prove that for small volumes, half-lines containing
$-\infty$ are minimizers. However, for $V=1/3$, we can consider a sequence of bounded intervals
of volume $V$ converging to $+\infty$ and whose perimeter tends to $2/9$.  A direct computation
shows that the half-line of volume $V$ containing $-\infty$ and any bounded interval of volume
$V$ have strictly greater perimeter. As consequence, there are no isoperimetric regions of
volume $1/3$ for this density.}
\end{example}

We say that a density $f$ is \emph{decreasing-increasing} if $(-f)$ is increasing-decreasing.
For these densities we have the following:

\begin{theorem}
\label{th:fnoinod}
Let $f$ be a decreasing-increasing density on $\rr$.  Then,
isoperimetric regions exist for any given volume and they are bounded
intervals in whose closure $f$ attains its minimum.
\end{theorem}

\begin{proof}
Take an open set $\Om$ with finite volume and a point $x_0$ where $f$ attains its minimum. It
is easy to check that the bounded interval containing $x_0$ and with the same volume as $\Om$
at both sides of $x_0$ has less perimeter than $\Om$. Finally, among intervals with fixed
volume containing $x_0$ in its closure there is one of least perimeter.
\end{proof}

Now, we give some applications and improvements of the previous results for the particular
cases of log-concave and log-convex densities.  We begin with the following corollary, which is
a direct consequence of Proposition~\ref{prop:fmonotona}, Theorem~\ref{th:fnodnoi}, and
elementary properties of concave functions.

\begin{corollary}
\label{th:logconcave}
Let $f$ be a log-concave density on $\rr$.  Then we have
\begin{itemize}
\item[(i)] If the total measure is finite, then minimizers of any
volume exist and they can be half-lines, unions of two disjoint half-lines, or bounded
intervals where the maximum of the density is attained.
\item[(ii)] If both ends have infinite measure,
then the density is constant and bounded intervals provide minimizers
of any given volume.
\item[(iii)] If the density has infinite volume but one end has finite
measure, then half-lines containing this end are the unique isoperimetric regions for fixed
volume.
\end{itemize}
\end{corollary}

Example \ref{re:houseroof} shows that all the different possibilities
in Corollary \ref{th:logconcave} (i) can appear.  C~Borell
(\cite[Corollary 2.2]{bor3}, see also \cite[Corollary 13.8]{bh})
proved that half-lines are always minimizers for a log-concave density
with finite total measure.  In the next corollary we give a different
proof of this fact showing also uniqueness of minimizers for strictly
log-concave densities.

\begin{corollary}
\label{cor:logconcave}
Let $f=e^\psi$ be a log-concave density on $\rr$ with finite total
measure.  Then, half-lines are always isoperimetric regions of any
given volume.  Moreover, if the density is strictly log-concave, then
half-lines are the unique minimizers.
\end{corollary}

\begin{proof}
We have to compare the perimeter of the candidates provided by Corollary~\ref{th:logconcave}
(i).  By taking complements we see that it is enough to compare the perimeter of bounded
intervals and half-lines of the same measure.  Fix an amount $V$ of volume.  Let $x_V$ be the
real number such that $\vol((x_V,+\infty))=V$.  For any $x\in (-\infty,x_V)$, let $y(x)>x$ be
the unique value satisfying $\vol((x,y(x)))=V$.  The perimeter of all bounded intervals
enclosing volume $V$ is given by the function $P(x)=f(x)+f(y(x))$. Clearly, $P(-\infty)$ and
$P(x_{V})$ represent the perimeter of the two half-lines of volume $V$.  As $y(x)$ is
increasing and the density is log-concave, we deduce that $P(x)$ is an absolutely continuous
function with left and right derivatives at every point.  In particular, the right derivative
$P'_r$ is given by
\[
P'_r(x)=f(x)\,\{\psi'_r(x)+\psi'_r(y(x))\},\qquad x\in (-\infty,x_V).
\]
On the other hand, as $\psi$ is concave, we get that $\psi'_r$ is non-increasing and hence
$P'_r(x)/f(x)$ is also non-increasing. Thus, $P(x)$ is monotonic or increasing--decreasing on
$(-\infty,x_V)$.  Anyway the infimum of $P(x)$ is achieved in a half-line of volume $V$.
Moreover, if $f$ is strictly log-concave, then the infimum of $P(x)$ is not attained on
$(-\infty,x_V)$, so that the half-lines are the unique minimizers.
\end{proof}

\begin{remark}
Two relevant examples in probability and statistics where Corollary \ref{cor:logconcave} is
applied are the standard Gaussian density $f(x)=e^{-\pi x^2}$ and the logistic density
$f(x)=e^{-x}\,(1+e^{-x})^{-2}$.  As indicated in \cite{b}, for these densities it is also
interesting to describe minimizers under a volume constraint of the functionals
$\vol(\Om+[-h,h])$ for any $h>0$.  In \cite[Remark 13.9]{bh} it is pointed out that half-lines
are solutions to this problem.  In higher dimension, we can consider the same problem with the
cube $[-h,h]^{n+1}$. It was shown in \cite[Theorem 1.1]{b} that half-spaces are minimizers for
any product measure $\mu^{n+1}$ in $\rrn$ provided $\mu$ is log-concave with finite total
volume (see also \cite[Corollary 15.3]{bh} for the particular case of the logistic density).
\end{remark}

Now we state a result similar to Corollary~\ref{th:logconcave} where we completely describe
isoperime\-tric regions for log-convex densities.

\begin{corollary}
\label{th:cpletedescriplogconv} Let $f=e^{\psi}$ be a log-convex density on $\rr$.  Then we
have
\begin{itemize}
\item[(i)] If both ends have infinite measure and
$f(-\infty)=f(+\infty)=+\infty$, then isoperimetric regions of any
volume exist and they are bounded intervals in whose closure the
density attains its minimum.  Moreover, if $f$ is strictly log-convex,
then we have uniqueness of minimizers for given volume.
\item[(ii)] If both ends have infinite measure but there is one end $E$ with
$f(E)<+\infty$, then the isoperimetric profile is constant and it is
approached or attained by a bounded interval going off to~$E$.
\item[(iii)] If one end has finite measure, then the half-lines con\-taining
this end are the unique minimizers for given volume.
\end{itemize}
\end{corollary}

\begin{proof}
The claim follows by using Proposition~\ref{prop:fmonotona} and Theorem~\ref{th:fnoinod}. The
uniqueness in statement (i) follows from the argument in the proof of Corollary
\ref{cor:logconcave} since strict convexity of the density implies that the perimeter of
bounded intervals with fixed volume achieves its minimum only at one point.
\end{proof}

As a direct consequence of the previous corollary and Example \ref{ex:real} we deduce the
follo\-wing result, which solves Conjecture~\ref{conj:conjecture} in dimension one.

\begin{corollary}
\label{cor:logconvsym} Let $f$ be a smooth, symmetric, strictly log-convex
density on $\rr$.  Then, for a given volume, the symmetric interval
of this volume is the unique minimizer.
\end{corollary}

The comparison arguments in this section allow to study the isoperimetric problem in
$[0,+\infty)$ or in a bounded interval $[a,b]$ with unimodal densities.  In these settings we
can prove similar results to Proposition~\ref{prop:fmonotona}, Theorem \ref{th:fnodnoi} and
Theorem~\ref{th:fnoinod}. The proofs are left to the reader.

\begin{theorem}
\label{th:idinhalfline} Let $f$ be a unimodal density on $[0,+\infty)$. Then we have
\begin{itemize}
\item[(i)] If $f$ is increasing, then the unique minimizers are the intervals $(0,x)$.
If $f$ is decreasing and $E=+\infty$ has finite measure, then the half-lines containing $E$ are
the unique minimizers. If $f$ is decreasing and $E$ has infinite measure, then the
isoperimetric profile equals $2f(E)$ and it is approached or attained by a bounded interval
going off to $+\infty$.
\item[(ii)] If $f$ is increasing-decreasing and a minimizer of given volume exists, then it
must coincide with an interval $(0,x)$, a half-line containing $+\infty$, a bounded interval
where $f$ attains its maximum or equals it one-sided minimum, or the complement of one of these
intervals. If it does not exist, the infimum perimeter is approached by a bounded interval
going off to $+\infty$.
\item[(iii)] If $f$ is decreasing-increasing, then minimizers of any measure exist and they are
bounded intervals in whose closure $f$ attains its minimum.
\end{itemize}
\end{theorem}

Now we shall state the corresponding result for the isoperimetric problem inside a bounded
interval $[a,b]$. Observe that in this case the existence of minimizers is ensured by
compactness.

\begin{theorem}
Let $f$ be a unimodal density on a bounded interval $[a,b]$. Then we have
\begin{itemize}
\item[(i)] If $f$ is monotonic, then any isoperimetric region is an interval
whose closure contains the boundary point of $[a,b]$ where $f$ attains its minimum.
\item[(ii)] If $f$ is increasing-decreasing, then a minimizer must coincide with an interval
whose closure contains a boundary point, an interval where $f$ attains its maximum or equals
its one-sided minimum, or the complement of one of these intervals.
\item[(iii)] If $f$ is decreasing-increasing, then any minimizer is an open interval in whose
closure $f$ attains its minimum.
\end{itemize}
\end{theorem}

The techniques in this section can also be applied to study the \emph{free boundary problem} in
$[0,+\infty)$ or $[a,b]$ which consists of finding global minimizers under a volume constraint
of the \emph{perimeter relative to} $(0,+\infty)$ or $(a,b)$, respectively. This means that the
boundary points of these intervals do not contribute to perimeter. For the case of
$[0,+\infty)$ we have:

\begin{theorem}
Let $f$ be a unimodal density on $[0,+\infty)$. Then we have
\begin{itemize}
\item[(i)] If $f$ is increasing, then the unique minimizers for the free boundary problem in
$[0,+\infty)$ are intervals of the form $(0,x)$. If $f$ is decreasing and $E=+\infty$ has
finite measure, then minimizers exist and they are half-lines containing $E$. If $E$ has
infinite measure and a minimizer exists, then it must coincide with an interval $(0,x)$ or a
bounded interval where $f$ equals its minimum. If a minimizer does not exist the infimum
perimeter equals $2f(E)$.
\item[(ii)] If $f$ is increasing-decreasing and a minimizer of given volume exists, then it is
an interval $(0,x)$, a half-line containing $+\infty$, a bounded interval where $f$ attains its
maximum or equals its right-side minimum, or the complement of one of these intervals. If it
does not exist, then the infimum perimeter is approached by a bounded interval going off to
$+\infty$.
\item[(iii)] If $f$ is decreasing-increasing then minimizers of any given volume are provided by
intervals $(0,x)$ or bounded intervals in whose closure $f$ attains its minimum.
\end{itemize}
\end{theorem}

For the free boundary problem in $[a,b]$, existence of minimizers is assured by compactness. As
to the description of isoperimetric regions, we can prove the following:

\begin{theorem}
Let $f$ be a unimodal density on $[a,b]$. Then
\begin{itemize}
\item[(i)] If $f$ is monotonic then the unique minimizers for the free boundary problem are
intervals whose closure contains a boundary point.
\item[(ii)] If $f$ is increasing-decreasing, then isoperimetric regions are provided by
intervals where $f$ attains its maximum, or whose closure contains a boundary point of $[a,b]$,
or complements of these intervals.
\item[(iii)] If $f$ is decreasing-increasing, then minimizers are intervals whose closure
contains a boundary point of $[a,b]$ or a value where the minimum of $f$ is attained, or
complements of these intervals.
\end{itemize}
\end{theorem}

\section{Isoperimetric inequality for the density $\exp(|x|^2)$}
\label{sec:logconvdensity}

In this last section of the paper we solve the isoperimetric problem in $\rrn$ with the radial
log-convex density $f(x)=\exp(c|x|^2)$, where $c$ is a positive constant.  Precisely, we will
prove that Conjecture \ref{conj:conjecture} holds for this density: round balls about the
origin provide isoperimetric regions of any given volume, like Euclidean space $(c=0)$ and
unlike Gauss space $(c<0)$.  As we pointed out in the Introduction, the proof combines Steiner
symmetrization in axis directions as was employed by L.~Bieberbach \cite{bieberbach} together
with Hsiang symmetrization \cite{hsiang}.  We will also show uniqueness by a detailed analysis
of the situation where an axis symmetrization of a minimizer produces a round ball.

The use of Steiner symmetrization in our setting is natural since the ambient density can be
seen as a rotationally invariant product measure.  Let us recall some facts about this
construction; see \cite[Section 3.2]{ros} for details.  Let $\Om$ be a compact set in $\rrn$.
Consider a hyperplane $\pi$ in $\rrn$ containing the origin. The restriction of the ambient
density to any straight line orthogonal to $\pi$ provides a smooth, symmetric, strictly
log-convex density. We define the \emph{symmetrization of} $\Om$ \emph{with respect to $\pi$}
as the set $\Om^*$ whose intersection with any straight line $R$ orthogonal to $\pi$ is the
isoperimetric region in $R$ of the same weighted length as $\Om\cap R$.  By Corollary
\ref{cor:logconvsym} this will be an interval centered at $\pi\cap R$.  It is clear that
$\Om^*$ is symmetric with respect to $\pi$.  The main property of this construction is that it
preserves volume (by Fubini's theorem) while decreasing perimeter.

\begin{lemma}[{\cite[Proposition 8]{ros}}]
\label{lem:symm}
For any hyperplane $\pi$ through the origin in $\rrn$, the Steiner
symmetrization $\Om^*$ of a compact set $\Om$ satisfies
$\vol(\Om^*)=\vol(\Om)$ and $P(\Om^*)\leq P(\Om)$.
\end{lemma}

Now, we will proceed to prove our main result in this section.

\begin{theorem}
\label{th:main}
In $\rrn$ with the density $f(x)=\exp(c|x|^2)$, $c>0$, round balls
about the origin uniquely minimize perimeter for given volume.
\end{theorem}

\begin{proof}
First observe that bounded minimizers of any given volume exist for this density by
Theorem~\ref{prop:exist1}.  Let us prove that round balls centered at the origin are
isoperimetric regions.  Take a minimizer $\Om$ of volume $V>0$.  We apply Steiner
symmetrization to $\overline{\Om}$ with respect to any coordinate hyperplane in $\rrn$ so that
we produce a minimizer $\Om^*$ which is symmetric with respect to any of these hyperplanes and
has connected boundary.  In particular, $\Om^*$ is centrally symmetric.  Thus any hyperplane
$\pi$ through the origin divides $\Om^*$ in two sets $\Om_i^*$ contained in the corresponding
open half-spaces $\pi_i$ and with the same volume. Note that the reflection with respect to
$\pi$ preserves the perimeter relative to any $\pi_i$ since the density is radial.  It follows
that $P(\Om_1,\pi_1)=P(\Om_2,\pi_2)$; otherwise, we would obtain by reflection a set with the
same volume as $\Om^*$ and strictly less perimeter.  Therefore each $\Om_i^*$ together with its
reflection is a new minimizer of volume $V$.  By the regularity properties in
Theorem~\ref{th:reg} and unique continuation for (real-analytic) generalized constant mean
curvature surfaces, $\ptl\Om^*$ is symmetric across any hyperplane $\pi$ through the origin. We
conclude that $\Om^*$ coincides with a ball centered at the origin.

To prove uniqueness, by induction it suffices to show that if symmetrization of a minimizer
$\Om$ with respect to a coordinate hyperplane $\pi$ produces a ball $B$, then $\Om$ is a ball.
We can suppose that $\pi=\{x_{n+1}=0\}$.  Let $D\sub\pi$ be the projection of $\Om$.  By
Theorem~\ref{th:reg} and Sard's theorem, for almost all $p\in D$ straight lines near $p$
orthogonal to $\pi$ intersect $\Sg=\ptl\Om$ transversally at a fixed even number of points
$p_i$, where $\Sg$ is the graph over $D_p\sub D$ of a smooth function $h_i$ (if we did not know
that $\Om$ is bounded, we would allow $p_i$ to be $\pm\infty$). Denote by $A\subeq D$ the set
of such points $p$. By the definition of Steiner symmetrization
\[
\sum_{i\,\,odd}\,\int_{h_i}^{h_{i+1}}f(x)\,dx=2\,\int_0^{h^*}f(x)\,dx\qquad\text{ on }\ D_p,
\]
where $h^*$ is the height function of $\ptl B$ with respect to $\pi$.
As a consequence
\[
\sum_{i\,\,odd}\,(f(h_{i+1})\,\nabla h_{i+1}-f(h_i)\,\nabla h_i)= 2f(h^*)\,\nabla
h^*\qquad\text{ on }\ D_p,
\]
so that we get
\[
\sum_j f(h_j)\,|\nabla h_j|\geq 2f(h^*)\,|\nabla h^*|\qquad\text{ on }\ D_p.
\]
On the other hand, by Corollary~\ref{cor:logconvsym} we have
\[
\sum_j f(h_j)\geq 2f(h^*)\qquad\text{ on }\ D_p,
\]
and equality holds if and only if the corresponding slice of $\Om$ is
an interval centered at $\pi$.

Now we apply Lemma~\ref{lem:technical} below with
$\alpha_j=f(h_j(p))$, $a_j=|\nabla h_j(p)|$, $\alpha=f(h^*(p))$ and
$a=|\nabla h^*(p)|$.  We get
\begin{equation}
\label{eq:th1}
\sum_j\, f(h_j(p))\,\sqrt{1+|\nabla h_j(p)|^2}\geq 2f(h^*(p))\,
\sqrt{1+|\nabla h^*(p)|^2},\qquad p\in A,
\end{equation}
with equality if and only if $|\nabla h_j(p)|=|\nabla h^*(p)|$ for any
$j$, and the slice of $\Om$ passing through $p$ is an interval
centered at $\pi$.

Finally we use the coarea formula and inequality \eqref{eq:th1} to obtain
\begin{align*}
P(\Om)=\int_{\Sg}f\,da&\geq\int_A\bigg(\sum_j\, f(h_j(p))\,
\sqrt{1+|\nabla h_j(p)|^2}\bigg)\,da
\\
&\geq\int_A 2f(h^*(p))\,\sqrt{1+|\nabla h^*(p)|^2}\,\,da=P(B),
\end{align*}
where in the last equality we have used that $D-A$ does not contribute to the perimeter of $B$.
As $\Om$ is a minimizer we have equality above and hence in \eqref{eq:th1} too.  It follows
that for every $p\in A$ the slice of $\Om$ passing through $p$ is a symmetric interval of the
same length as the corresponding slice for $B$.  Thus, up to a set of measure zero, $\Om$
coincides with a round ball about the origin.
\end{proof}

\begin{lemma}
\label{lem:technical}
Suppose that we have finitely many nonnegative real numbers with
$\sum_{j}\alpha_{j}\,a_{j}\geq 2\alpha\,a$ and $\sum_{j}\alpha_{j}\geq
2\,\alpha$.  Then the following inequality holds
\[
\sum_{j}\alpha_{j}\,\sqrt{1+a_{j}^2}\geq 2\alpha\sqrt{1+a^2},
\]
with equality if and only if $a_{j}=a$ for every $j$ and
$\sum_j\alpha_{j}=2\alpha$.
\end{lemma}

\begin{proof}
The function $g(x)=\sqrt{1+x^2}$ is strictly convex with\, $0<g'(x)<1$
for any $x>0$.  Let $\alpha_0=\sum_{j}\alpha_j$.  We claim that
\[
\sum_{j}\frac{\alpha_{j}}{\alpha_{0}}\,\,g(a_{j})\geq
g\big(\sum_{j}\frac{\alpha_{j}}{\alpha_{0}}\,\,a_j\big)\geq
g\big(\frac{2\alpha}{\alpha_0}\,a\big)\geq \frac{2\alpha}{\alpha_0}\,\, g(a).
\]
The first inequality holds because $g$ is convex.  The second and third inequalities come from
the fact that $0<g'(x)<1$ for $x>0$.  If equality holds in the second inequality, then
$\sum_{j}\alpha_{j}\,a_{j}=2\alpha\,a$.  If equality holds in the third inequality too, then
$2\alpha=\alpha_0=\sum_{j}\alpha_j$. If equality holds in the first inequality as well, then
$a_{j}=a$ for every $j$.
\end{proof}

We finish the paper with an eigenvalues comparison theorem obtained as
a consequence of the isoperimetric inequality in
Theorem~\ref{th:main}.  For a smooth bounded domain $\Om$ in $\rrn$,
we consider the second order differential operator $L$ on
$C^\infty_{0}(\Om)$ whose invariant measure has density
$f(x)=\exp(c|x|^2)$ $(c\geq 0)$, namely
\begin{equation}
\label{eq:operator}
(Lu)(x)=(\Delta u)(x)-2c\escpr{x,(\nabla u)(x)},\qquad u\in
C^\infty_{0}(\Om),\quad x\in\Om,
\end{equation}
where $\Delta$ denotes the Euclidean Laplace operator on $\Om$.

\begin{corollary}
\label{cor:eigenvalues}
Let $\Om$ be a smooth bounded domain in $\rrn$.  Then, the lowest
non-zero eigenvalue $\lambda_{1}(\Om)$ for the second order
differential operator \eqref{eq:operator} with Dirichlet boundary
condition on $\ptl\Om$ satisfies
\[
\lambda_{1}(\Om)\geq\lambda_{1}(B),
\]
where $B$ is the round ball centered at the origin with the same
volume as $\Om$ for the density $f(x)=\exp(c|x|^2)$, $c>0$.
Moreover, equality holds if and only if $\Om=B$.
\end{corollary}

\begin{proof}
The comparison is an adaptation of the symmetrization technique used
to prove the Faber-Krahn Inequality in $\rrn$ (see \cite[p.
87]{chavel2}), which corresponds to the desired inequality for the
case $c=0$.
\end{proof}

\providecommand{\bysame}{\leavevmode\hbox to3em{\hrulefill}\thinspace}
\providecommand{\MR}{\relax\ifhmode\unskip\space\fi MR }
\providecommand{\MRhref}[2]{%
   \href{http://www.ams.org/mathscinet-getitem?mr=#1}{#2}
}
\providecommand{\href}[2]{#2}

\end{document}